\documentclass[11pt, oneside]{article}   	% use "amsart" instead of "article" for AMSLaTeX format
\usepackage{geometry}                		% See geometry.pdf to learn the layout options. There are lots.
\usepackage{graphicx}				% Use pdf, png, jpg, or eps§ with pdflatex; use eps in DVI mode
\usepackage{amssymb}

\newtheorem{theorem}{Theorem}

\newtheorem{remk}[theorem]{Remark}

\def\FullBox{\hbox{\vrule width 8pt height 8pt depth 0pt}}

\def\qed{\ifmmode\qquad\FullBox\else{\unskip\nobreak\hfil
\penalty50\hskip1em\null\nobreak\hfil\FullBox
\parfillskip=0pt\finalhyphendemerits=0\endgraf}\fi}

\def\qedsketch{\ifmmode\Box\else{\unskip\nobreak\hfil
\penalty50\hskip1em\null\nobreak\hfil$\Box$
\parfillskip=0pt\finalhyphendemerits=0\endgraf}\fi}

\title{A  Cameron and Erd\H os conjecture on counting primitive sets}
\author{Rodrigo Angelo - rangelo@princeton.edu}
\begin{document}

\maketitle

%\section{}
%\subsection{}

\abstract{Let $f(n)$ count the number of  subsets of $\{1,...,n\}$ without an element dividing another. In this paper I show that $f(n)$ grows like the $n$-th power of some real number, in the sense that $\lim_{n\rightarrow \infty}f(n)^{1/n}$ exists. This confirms a conjecture of Cameron and Erd\H os, proposed in a paper where they studied a number of similar problems, including the well known "Cameron-Erd\H os Conjecture" on counting sum-free subsets.
\paragraph{}
\paragraph{}

Let $f(n)$ be the number of subsets of $[n]=\{1,...,n\}$ such that no element divides another (call these sets primitive). One easily notices $2^n \ge f(n)\ge 2^{n/2}$ (since subsets of the second half are all primtive), motivating Cameron and Erd\H os to question whether there is an exact real number characterizing the exponential growth of this function [1]. We confirm their conjecture:

\paragraph{Theorem:} $\lim_{n\rightarrow \infty}f(n)^{1/n}$ exists

\paragraph{Proof:}  We will study the auxiliary and more structured $f(n,k)$, which we define to be the number of subsets of $[n]$ such that no two elements have an integer ratio for which all prime factors are at most $p_k$ (the $k$-th prime number). Call these sets $k$-core. The crux of the proof will be a little argument that shows that if for each $k$ $\lim_{n\rightarrow \infty} f(n,k)^{1/n} $ exists, then $\lim_{n\rightarrow \infty}f(n)^{1/n}$ also exists. This is somewhat surprising, because if one doesn't think about this the right way it may seem that it is necessary to send $k$ to infinity together with $n$ in order to obtain the desired limit. That said, we divide the proof into two parts:

Part 1: If we assume that for each $k$ $\lim_{n\rightarrow \infty} f(n,k)^{1/n} = \alpha_k$ exists, then the $\alpha_k$ decrease to some limit $\alpha$ and $\lim_{n\rightarrow \infty}f(n)^{1/n}$ exists and is equal to $\alpha$

Part 2: For each fixed $k$, $\lim_{n\rightarrow \infty}$ $f(n,k)^{1/n} $ in fact exists

\paragraph{Proof of part 1:}

Let $\lim_{n\rightarrow \infty} f(n,k)^{1/n} = \alpha_k$. Clearly $f(n,k+1) \le f(n,k)$ (because the condition of being $k+1$-core is more restrictive than the condition of being $k$-core). By taking $1/n$ powers and limits we get that $\alpha_k$ are a decreasing sequence. Since they are non-negative it follows that the $\alpha_k$ must have a limit $\alpha$. 

Since $f(n) \le f(n,k)$ for all $k$, we get $\limsup f(n)^{1/n} \le \alpha_k$ for each $k$, which gives $\limsup f(n)^{1/n} \le \alpha$. 

Now we need an inequality for the other side.

For that, we notice that for a $k$-core subset of $[n]$, if the elements less than $\frac{n}{k}$ are removed we get a primitive subset (because the ratio of any two remaining elements is less than $k$, so if one element divided another their ratio would be an integer less than $k$, which has all prime factors less than $p_k$, which contradicts that the original set is $k$-core). So this operation maps $k$-core sets to primitive sets. Also, it is clear that this operation maps at most $2^{n/k}$ sets to the same set (because two sets mapped to the same set may disagree only on the first $n/k$ elements). This gives the inequality:

$$f(n,k) \le 2^{n/k}f(n)$$

By taking $1/n$ power and taking $n$ to infinity this gives:

$$\liminf f(n)^{1/n} \ge \alpha_k 2^{-1/k},$$

for all $k$. By making $k \rightarrow \infty$ we get $\liminf f(n)^{1/n} \ge \alpha$. So

$$\alpha \le \liminf f(n)^{1/n} \le \limsup f(n)^{1/n} \le \alpha,$$

which completes the proof that $\lim_{n \rightarrow \infty} f(n)^{1/n}$ exists and is equal to $\alpha$.

\paragraph{Proof of part 2:}

Fix $k$. Let $S = \{p_1,...,p_k\}$ and $D = p_1...p_k$ be the product of the first $k$ primes. Each integer can be written uniquely as a product $aR$ where $a$ only has prime factors in $S$ and $(R,D)=1$. Integers with distinct values of $R$ cannot have an integer ratio with prime factors in $S$. So we partition the integers in $[n]$ according to their value of $R$, and the total number of $k$-core subsets of $[n]$ is just the product of the number of $k$-core subsets of each part. We also notice that each part consists of the integers of the form $aR$, where a runs over the integers less than $\frac{n}{R}$ with all prime factors in $S$. Hence if we define $P_k(x)$ to be the number of $k$-core (or simply primitive) subsets of the set of integers $\le x$ with all prime factors in $S$ we get:

$$f(n,k) = \prod_{1\le R \le n, (R,D)=1} P_k(\lfloor \frac{n}{R} \rfloor).$$

Now set $\epsilon>0$ to be chosen later. We first want to show that the first $\epsilon n$ terms of this product do not contribute substantially. For these terms we use the bound:

$$P_k(x) \le 2^{(1+\log x)^k}$$

(we obtain this by bound $P_k(x)$ above by the number of subsets of the set of integers less than $x$ with all prime factors in S, and we bound the size of this set by $(1+\log x)^k$ by noticing that each $p_1^{a_1}...p_k^{a_k} \le x$ is associated to a distinct $k$-tuple $(a_1,...a_k)$ with $a_i \le 1+\log x$). Hence:

$$\prod_{1\le R \le \epsilon n, (R,D)=1} P_k(\lfloor \frac{n}{R} \rfloor) \le  \prod_{1\le R \le \epsilon n} 2^{(1+\log \frac{n}{R})^k} \le 2^{\epsilon n (1+\log n)^k}.$$

The product of the first $\epsilon n$ terms is also $\ge 1$, so we get:

$$f(n,k) = 2^{O(\epsilon n (1+\log n)^k)} \prod_{\epsilon n\le R \le n, (R,D)=1} P_k(\lfloor \frac{n}{R} \rfloor).$$

Now $\frac{n}{R}$ is always between 1 and $\frac{1}{\epsilon}$. For each integer $l$ between 1 and $\frac{1}{\epsilon}$ there are $n(\frac{1}{l}-\frac{1}{l+1}) + O(1)$ integers $R$ from $\epsilon n$ to $n$ with $\lfloor \frac{n}{R} \rfloor=l$. And this is a run of consecutive numbers, so $n(\frac{1}{l}-\frac{1}{l+1}) \frac{\phi(D)}{D}+O(D)$ of these numbers are prime with $D$ ($\phi(D)$ is the Euler function). Hence:

$$f(n,k)^{1/n} = 2^{O(\epsilon  (1+\log n)^k)} \prod_{1\le l \le \frac{1}{\epsilon}} P_k(l)^{(\frac{1}{l}-\frac{1}{l+1}) \frac{\phi(D)}{D}+\frac{O(D)}{n}}$$  

$$= 2^{O(\epsilon  (1+\log n)^k+\frac{D (1+\log 1/\epsilon)^k}{\epsilon n})}\prod_{1\le l \le \frac{1}{\epsilon}} P_k(l)^{(\frac{1}{l}-\frac{1}{l+1}) \frac{\phi(D)}{D}}.$$

Here we used the bound $P_k(l) \le 2^{(1+\log l)^k}\le 2^{(1+\log 1/\epsilon)^k}$ again. Finally we choose $\epsilon = \frac{1}{\sqrt{n}}$. By making $n \rightarrow \infty$ both the error terms go to zero, and the number of terms in the product go to infinity, so in order to prove $\lim_{n\rightarrow \infty} f(n,k)^{1/n}$ exists it is enough to show:

$$\prod_{l=1}^{\infty} P_k(l)^{(\frac{1}{l} - \frac{1}{l+1})\frac{\phi(D)}{D}}$$

is a convergent product (and the limit will be equal to this product). Indeed, by the same bound for $P_k(x)$, it is enough to prove $\sum_{l=1}^{\infty} \frac{(1+\log l)^k}{l(l+1)}$ is convergent, which is true. Hence the proof is complete.

\paragraph{}

Unfortunately my attempts up to now have failed to find the value of $\alpha$ (in some reasonable sense). What seems to happen is that this solution essentially reduces the limit to a "somoothed" version of the limit in terms of the $P_k$ which is guaranteed to converge - but because we don't know much else about the $P_k$ attempts to find the limit end up circular. It is also amusing to notice that if one looks only at the infinite product formula we found for $\alpha_k$ it is not obvious that these form a decreasing sequence - one needs the "combinatorial" argument from part 1 to stablish that, and this seems to be a considerable barrier to making sense out of the limit of $\alpha_k$ trough this formula.

\paragraph{References}
\paragraph{}
[1] Cameron, P.J; Erd\H os, P. \textit{On the number of sets of integers with various properties},
Number theory (Banff, AB, 1988), 61-79, de Gruyter, Berlin, 1990.
\end{document}